\documentclass[letterpaper,10pt]{amsart}

\usepackage{epsfig}
\usepackage{amssymb}
\usepackage{amsfonts}
\usepackage{amsmath}
\usepackage{graphicx}
\usepackage{mathrsfs}
\usepackage{stmaryrd}
\usepackage{amsthm}
\usepackage{textcomp}
\usepackage[all]{xy}
\usepackage[top=1.25in, bottom=1.25in, left=1.25in, right=1.25in]{geometry}
\usepackage{xcolor}
\usepackage[colorlinks=true]{hyperref}
\hypersetup{urlcolor=blue, citecolor=red, linkcolor=blue}
\usepackage[colorinlistoftodos]{todonotes}
\usepackage{color}

\newtheorem{theorem}{Theorem}

\newtheorem{corollary}{Corollary}

\def \pa{\partial}

\newcommand{\ud}{\mathrm{d}}

\title[Escobar-Yamabe compactifications for Poincar\'{e}-Einstein manifolds]{Escobar-Yamabe compactifications for Poincar\'{e}-Einstein manifolds and rigidity theorems}
\author{Xuezhang Chen}
\address{Department of Mathematics \& IMS, Nanjing University, Nanjing
210093, P. R. China}
\email{xuezhangchen@nju.edu.cn}
\author{Mijia Lai}
\address{School of Mathematical Sciences, Shanghai Jiao Tong University, Shanghai 200240, P. R. China}
\email{laimijia@sjtu.edu.cn}
\author{Fang Wang}
\address{School of Mathematical Sciences, Shanghai Jiao Tong University, Shanghai 200240, P. R. China}
\email{fangwang1984@sjtu.edu.cn}

\thanks{Chen's research was supported by NSFC (No.11771204), A Foundation for the Author of National Excellent Doctoral Dissertation of China (No.201417) and start-up grant of 2016 Deng Feng plan B at Nanjing University. Lai's research was supported in part by National Natural Science Foundation of China No.11501360.
 Wang's research was supported in part by National Natural Science Foundation of China No. 11401377. }
\date{}
\begin{document}

\begin{abstract}
Let $(X^{n},g_+) $ $(n\geq 3)$ be a Poincar\'{e}-Einstein manifold which is $C^{3,\alpha}$ conformally compact with conformal infinity $(\partial X, [\hat{g}])$. On the conformal compactification $(\overline{X}, \bar g=\rho^2g_+)$ via some boundary defining function $\rho$, there are two types of Yamabe constants: $Y(\overline{X},\partial X,[\bar g])$ and $Q(\overline{X},\partial X,[\bar g])$. (See definitions (\ref{def.type1}) and (\ref{def.type2})). In \cite{GH},  Gursky and Han gave an inequality between $Y(\overline{X},\partial X,[\bar g])$ and $Y(\partial X,[\hat{g}])$. In this paper, we first show that the equality holds in Gursky-Han's theorem if and only if $(X^{n},g_+)$ is isometric to the standard hyperbolic space $(\mathbb{H}^{n}, g_{\mathbb{H}})$. Secondly, we derive an inequality between $Q(\overline{X},\partial X,[\bar g])$ and $Y(\partial X, [\hat g])$, and show that the equality holds if and only if $(X^{n},g_+)$ is isometric to $(\mathbb{H}^{n}, g_{\mathbb{H}})$. Based on this, we give a simple proof of the rigidity theorem for Poincar\'{e}-Einstein manifolds with conformal infinity being conformally equivalent to the standard sphere.
\end{abstract}

\maketitle

\section{Introduction}
Let $\overline{X}$ be an $n$-dimensional smooth connected compact manifold with boundary.
Denote by $X$ the interior of  $\overline{X}$ and by $\partial X$ the boundary of $\overline{X}$.
Let $g_+$ be a Remannian metric in $X$. 
We call
$(X, g_+)$ a $C^{3,\alpha}$ conformally compact Poincar\'{e}-Einstein manifold if $g_{+}$ is a complete metric in $X$ satisfying 
$$
\mathrm{Ric}_{g_+}=-(n-1)g_+,
$$
and $\bar g=\rho^2 g_+$ can be $C^{3,\alpha}$ extended to $\overline{X}$ by some boundary defining function $\rho$.
Here $\rho$ is said to be a smooth boundary defining function if
$$
\textrm{$\rho\in C^{\infty}(\overline{X})$,  $\rho>0$ in X, $\rho=0$ on $\partial X$ and $\ud \rho\neq 0$ on $\partial X$. }
$$
Thus $(X,g_{+})$ is conformally
compactified to  $(\overline{X}, \bar{g})$, which is a compact Riemannian manifold with boundary. It is clear that $\bar g$ and the induced metric $\hat{g}=\bar g|_{\partial X}$ depend on the choice of $\rho$, but their conformal classes do not. In particular, we call
 $(\partial X, [\hat{g}])$ the conformal infinity of $(X, g_+)$.

Our work in this paper is motivated by the recent work of Gursky-Han~\cite{GH}. To describe their and our
results explicitly, we first recall the Yamabe constant on the conformal infinity. It is defined as follows:
$$
Y(\partial X,[\hat{g}]):=\inf_{\hat{h} \in[\hat{g}]}
\frac{ \int_{\partial X} R_{\hat{h}}   \ud S_{\hat{h}}    }
{(\int_{\partial X} \ud S_{\hat{h}} )^{\frac{n-3}{n-1}} },
$$
where $R_{\hat{h}}$ is the scalar curvature of metric $\hat{h}$ on $\partial X$.
This is related to the classical Yamabe problem on closed manifolds, which was completely solved.
See \cite{Au1, Tru, Sch, LP} and many other works.
For a compact manifold with boundary $(\overline{X}, \partial X, [\bar g])$, there are two types of Yamabe constants defined as follows:
\begin{equation}\label{def.type1}
Y(\overline{X},\partial X,[\bar  g]):=\inf_{\bar{h} \in[\bar g]}\frac{\int_X R_{\bar{h} }\ud V_{\bar{h} }
+2\int_{\partial X}H_{\bar{h} }\ud S_{\bar{h} }}
{(\int_{ X} \ud V_{ \bar{h} })^{\frac{n-2}{n}}};
\end{equation}
\begin{equation}\label{def.type2}
Q(\overline{X},\partial X,[\bar g]):=\inf_{\bar{h} \in[\bar g]}\frac{\int_X R_{\bar{h} }\ud V_{\bar{h} }+2\int_{\partial X}H_{\bar{h} }\ud S_{\bar{h} }}{(\int_{\pa X} \ud S_{\bar{h} })^{\frac{n-2}{n-1}}}.
\end{equation}
where $R_{\bar{h} }$ is the scalar curvature of metric $\bar{h} $ and $H_{\bar{h} }$ is the mean curvature of the boundary.
The related Yamabe problem for compact manifolds with boundary was intensively studied in the past half century.

For $Y(\overline{X},\partial X,[\bar  g])$,  Escobar showed in \cite{Es1} that for $n\geq 3$
\begin{equation}\label{ineq.Y}
Y(\overline{X},\partial X,[\bar  g])\leq Y(S^n_+, \partial S^n_+, [g_{S^n}]) =n(n-1)\left(\frac{1}{2}\omega_n\right)^{\frac{2}{n}},
\end{equation}
where $S^n_+$ is the upper hemisphere, $g_{S^n}$ is the round metric and $\omega_n$ is the volume of unit $n$-sphere. When the inequality is strict, then $Y(\overline{X},\partial X,[\bar  g])$ is attained by a metric $\bar{g}$ of constant scalar curvature with minimal boundary, i.e.,
$$
\begin{aligned}
& R_{\bar{g}}= Y(\overline{X},\partial X,[\bar  g]) \mathrm{Vol}(\overline{X}, \bar{g})^{-\frac{2}{n}}, \quad &\mathrm{in}\ X,
\\
& H_{\bar{g}}=0, \quad &\mathrm{on}\ \partial X.
\end{aligned}
$$

Escobar also showed that the inequality in (\ref{ineq.Y}) is strict if $(\overline{X},\partial X, \bar{g})$ is not conformally equivalent to $(S^n_+, \partial S^n_+, g_{S^n})$ and  $3\leq n\leq 5$, or $n\geq 6$ and $\partial X$ is not umbilic. In \cite{BCh}, Brendle-Chen considered the remaining case: $n\geq 6$ and $\partial X$ is umbilic. They verified the remaining case subject to the validity of the Positive Mass Theorem (PMT). By \cite{Wi}, the PMT is valid for spin manifold if $n\geq 6$. Sumarizing their results and applying them to  $C^{3,\alpha}$ conformal compactification of Poincar\'{e}-Einstein manifolds, which always have umbilic boundary, we get the following:

\begin{theorem}\label{thm.Y}
Let $(X^n, g_+)$ be a $C^{3,\alpha}$ conformally compact Poincar\'{e}-Einstein manifold satisfying one of the followings:
\begin{itemize}
\item[(a)] the dimension $3\leq n\leq 5$;
\item[(b)] the dimension $n\geq 6$ and $X$ is spin.
\end{itemize}
Then there is a conformal compactification $\bar{g}=\rho^2g_+$ satisfying
\begin{itemize}
\item[(1)] the scalar curvature $R_{\bar{g}}$ is constant;
\item[(2)] the boundary is totally geodesic, which implies that $H_{\bar{g}}=0$;
\end{itemize}
moreover the Yamabe constant $Y(\overline{X},\partial X, [\bar{g}])$ is achieved.
We call such  $(\overline{X}, \partial X, \bar{g})$ the  \textbf{first type Escobar-Yamabe compactification}  of $(X^n, g_+)$.
\end{theorem}

For $Q(\overline{X},\partial X,[\bar  g])$, in \cite{Es2} Escobar proved that for  $n\geq 3$, there holds
\begin{equation}\label{ineq.Q}
Q(\overline{X},\pa X, [\bar g]) \leq Q(\overline{B^n},S^{n-1},[g_{\mathbb{R}^n}]),
\end{equation}
where $B^n$ is the unit Euclidean ball and $g_{\mathbb{R}^n}$ is the Euclidean metric.
Moreover, when the inequality is strict, $Q(\overline{X},\pa X, [\bar g])$ is achieved by a scalar flat metric $\bar{g}$ with boundary having constant mean curvature, i.e.,
$$
\begin{aligned}
&R_{\bar{g}}=0, \quad& \mathrm{in}\ \overline{X},
\\
&H_{\bar{g}}=\frac{1}{2}Q(\overline{X},\partial X,\bar g)\mathrm{Vol}(\pa X,\hat g)^{-\frac{1}{n-1}}, \quad& \mathrm{on}\ \partial X,
\end{aligned}
$$
where $\hat{g}=\bar{g}|_{\partial X}$.

Assuming that $(\overline{X},\pa X,[\bar g])$ is not conformally equivalent to $(\overline{B^n},S^{n-1},[g_{\mathbb{R}^n}])$ and $Q(\overline{X},\partial X,\bar g)> -\infty$, Escobar was able to verify that the inequality \eqref{ineq.Q} is strict if  $n=3$; or $n=4,5$ and $\pa X$ is umbilic; or $n\geq 6$ and $X$ is locally conformally flat with umbilic boundary; or $n\geq 6$ and $\pa X$ has a non-umbilic point (in this case, $n=6$ was proved in \cite{Es3}). When $n=4,5$ and $\pa X$ has a non-umbilic point, the strict inequality was verified by Marques \cite{Mar2}. In \cite{Ch}, S. Chen considered the remaining case of $n\geq 6$ and $\pa X$ is umbilic. Moreover, Marques \cite{Mar1} developed an important tool of conformal Fermi coordinates, which plays the same role as conformal normal coordinates in the Yamabe problem. He then constructed appropriate test functions without using PMT in the umbilic boundary case: (i) $n=8$ and the Weyl tensor of $\pa X$ is nonzero at some boundary point; (ii) $n\geq 9$ and the Weyl tensor of $\overline X$ is nonzero at some boundary point.
See also \cite{Al} for a flow approach.  The work of S. Chen and that of Almaraz are particularly relevant to our setting. Applying their results, we have 
\begin{theorem}\label{thm.Q}
Let $(X^n,g_+)$ be $C^{3,\alpha}$ conformally compact Poincar\'e-Einstein metric satisfying one of the following:
\begin{enumerate}
\item[(a)] the dimension $3 \leq n\leq 7$;
\item[(b)] the dimension $n\geq 8$ and $X$ is spin;
\item[(c)] the dimension $n\geq 8$ and $X$ is locally conformally flat.
\end{enumerate}
Then there exists a conformal compactification $\bar g=\rho^2 g_+$ satisfying
\begin{itemize}
\item[(1)] the scalar curvature $R_{\bar{g}}=0$;
\item[(2)] the boundary mean curvature $H_{\bar g}$ is a constant,
\end{itemize}
and the Yamabe constant $Q(\overline{X},\partial X, [\bar{g}])$ is achieved.
We call such $(\overline{X}, \partial X, \bar{g})$ the \textbf{second type Escobar-Yamabe compactification}  of $(X^n, g_+)$.
\end{theorem}

Readers are referred to \cite{ChS} for more details  and generalizations on the Yamabe problem for compact manifolds with boundary.

In \cite{GH}, the authors derived an inequality between $Y(\overline{X},\partial X,[\bar g])$ and $Y(\partial X,[\hat{g}])$ using the first type Escobar-Yamabe compactification of $(X, g_+)$. The inequality involves an isoperimetric ratio $I(\overline{X}, \partial X, \bar{g})={\mathrm{Vol}(\partial {X}, \hat{g})^n} / {\mathrm{Vol}(\overline{X}, \bar{g})^{n-1}}$.  An interesting application of this inequality is the nonexistence of Poincar\'{e}-Einstein metric in the unit Euclidean ball $B^8$ with conformal infinity being $(S^7,[g])$ for infinitely many conformal classes $[g]$.  Their inequality is the following:
\begin{theorem}[Gursky-Han]\label{thm.GH}
Let $(X^n, g_+)$ be a Poincar\'{e}-Einstein manifold satisfying the hypotheses of Theorem \ref{thm.Y}. Let $(\overline{X}, \partial X, \bar{g})$ be the first type Escobar-Yamabe compactification and $\hat{g}=\bar{g}|_{\partial X}$.
\begin{itemize}
\item If the dimension $n\geq 4$, then
$$
Y(\overline{X},\partial X,[\bar g]) \cdot I(\overline{X}, \partial X, \bar{g})^{\frac{2}{n(n-1)}}
\geq \frac{n}{n-2} Y(\partial X, [\hat{g}]);
$$
\item If the dimension $n=3$, then
$$
Y(\overline{X},\partial X,[\bar g])\cdot I(\overline{X}, \partial X, \bar{g})^{\frac{1}{3}}\geq 12\pi \chi(\partial X).
$$
\end{itemize}
The equality is sharp and achieved when $(\overline{X}, \partial X, \bar{g})$ is the standard hemisphere. Moreover,
if the equality holds then $\bar{g}$ is Einstein and $\hat{g}$ has constant scalar curvature.
\end{theorem}

In this paper, we first show that the equality case in Theorem \ref{thm.GH} actually yields a rigidity theorem. We prove the following

\begin{theorem}\label{thm.2}
Let $(X^n, g_+)$ be a Poincar\'{e}-Einstein manifold satisfying the hypotheses of Theorem \ref{thm.Y}, with first type Escobar-Yamabe compactification $(\overline{X}, \partial X, \bar{g})$. Let $\hat{g}=\bar{g}|_{\partial X}$.
If the equality in Theorem \ref{thm.GH} holds, then $(X^n, g_+)$ is isometric to the standard hyperbolic space $(\mathbb{H}^{n}, g_{\mathbb{H}})$.
\end{theorem}

Secondly, we study the relation between $Q(\overline{X},\partial X,[\bar g])$ and $Y(\partial X,[\hat{g}])$. By working on the second type Escobar-Yamabe compactification, we derive an inequality only involving $Q(\overline{X},\partial X,[\bar g])$ and $Y(\partial X,[\hat g])$. The difference from Gursky-Han's inequality is that no volume ratio is involved and hence the inequality is conformally invariant.
Moreover, we show that the equality holds if and only if $(X, g_+)$ is isometric to the hyperbolic space.
\begin{theorem}\label{thm.1}
Let $(X,g_+)$ be  Poincar\'{e}-Einstein metric which satisfies hypotheses of Theorem \ref{thm.Q}.
Let $\bar{g}=\rho^2g_+$ be any conformal compactification and $\hat{g}=\bar{g}|_{\pa X}$ be the induced metric on  $\partial X$.
\begin{itemize}
\item If the dimension $n\geq 4$, then
\begin{equation}\label{ineq:Q_Y_n>3}
Y(\pa X, [\hat g])\leq \frac{n-2}{4(n-1)} Q(\overline{X},\pa X,[\bar g])^2;
\end{equation}
\item If the dimension $n=3$,  then
\begin{equation}\label{ineq:Q_Y_n=3}
32\pi \chi (\pa X)\leq Q(\overline{X},\pa X,[\bar g])^2.
\end{equation}
\end{itemize}%
In both cases, the equality holds if and only if $(X, g_+)$ is isometric to the hyperbolic space $(\mathbb{H}^n,g_{\mathbb{H}})$.
\end{theorem}

Similar to Gursky-Han's result, the assumptions in Theorem \ref{thm.Q} and Theorem \ref{thm.1} are mainly because of the limited validity of positive mass theorem (PMT).
In the proof of Theorem \ref{thm.1} we need to choose the second type Yamabe compactification, i.e. the Escobar-Yamabe metric of zero scalar curvature and constant boundary mean curvature such that $Q(\overline{X},\partial X,[\bar g])$ is achieved. The existence of such metric relies on PMT at the present stage. If the higher dimensional PMT is valid, then this assumption is expected to be removed. 

An application of Theorem \ref{thm.1} is to provide a simple proof of the rigidity theorem for Poincar\'{e}-Einstein manifolds.

\begin{corollary}\label{cor1}
Let $(X,g_+)$ satisfy the same hypotheses as in Theoerem \ref{thm.Q}.
If the conformal infinity of $(X,g_+)$ is conformally equivalent to the standard sphere $(S^{n-1},g_{S^{n-1}})$,  then $(X,g_+)$  is isometric to the hyperbolic space $(\mathbb{H}^n,g_{\mathbb{H}})$.
\end{corollary}

This rigidity theorem was proved in \cite[Theorem 1.1]{Qi1} for $4 \leq n\leq 7$, eventually it was extended to any dimension in \cite[Corollary 1.5]{LQS}, see also \cite{ST,DJ}. Here as a corollary of Theorem \ref{thm.1}, the proof is much shorter and more straight-forward. The disadvantage here is that we need a pre-assumption as in Theorem \ref{thm.1} such that the PMT in the boundary case works.

Here we would like to mention that a very closely related project to the second Escobar-Yamabe compactification is an isoperimetric inequality over scalar flat conformal class, which was studied by Hang-Wang-Yan \cite{HWY1,HWY2}. They proposed a conjecture therein: Let $(\overline {X^n},\pa X, \bar g)$ be a smooth compact Riemannian manifold of dimension $n \geq 3$ with smooth boundary, if $(\overline {X^n}, \pa X, [\bar g])$ is not conformally diffeomorphic to $(\overline{B^n},S^{n-1},[g_{\mathbb{R}^n}])$, then $\Theta(\overline X,\pa X, [\bar g])>\Theta (\overline{B^n},S^{n-1},[g_{\mathbb{R}^n}]),$ where
$$
\Theta(\overline X,\pa X, [\bar g])=
\sup_{\tilde g \in [\bar g], R_{\tilde g}=0} I(\overline{X}, \partial X, \tilde{g})^{\frac{-1}{n(n-1)}}.
$$
Furthermore, they proved that assuming this conjecture is true, if the first Dirichlet eigenvalue of conformal Laplacian of metric $\bar g$ is positive, then $\Theta(\overline X,\pa X, [\bar g])$ is achieved by some conformal metric with zero scalar curvature. Readers are referred to \cite{JX} for very recent progress in the above conjecture.

This paper is organized as follows. In Section \ref{sec2}, we recall Graham's work for singular Yamabe metrics, as well as an asymptotic calculus for boundary defining function $\rho$. In Section \ref{sec3}, we prove Theorem \ref{thm.2} by a flow analysis on the first type Escobar-Yamabe compactification and a generalization of Obata's rigidity theorem. In Section \ref{sec4}, we study the geometric formulae for the second type Escobar-Yamabe compactification and prove Theorem \ref{thm.1}.
In Section \ref{sec5}, we prove Corollary \ref{cor1} as a simple application of Theorem \ref{thm.1}.

\textbf{Acknowledgement}: The authors want to thank professor Paul Yang for valuable comments on this work.

\vspace{0.2in}
\section{Singular Yamabe and Asymptotics}\label{sec2}

We recall in this section some asymptotic computations by Graham \cite{Gr1} as well as the equations for scalar curvatures and Ricci curvatures under conformal change.

Let $(X^n, g_+)$ be a $C^{3,\alpha}$  conformally compact Poincar\'{e}-Einstein and $\bar{g}=\rho^2 g_+$ is a compactification. Let $\hat{g}=\bar{g}|_{\partial X}$. The normal exponential map
$\mathrm{exp}:[0,\delta)_r\times \partial X\longrightarrow \overline{X}$ relative to $\bar{g}$ is a diffeomorphism onto a neighbourhood of $\partial X$, with respect to which $\bar{g}$ takes the form
$$
\bar{g} =  \ud r^2+h_r,
$$
where $r(x)=d_{\bar g}(x,\pa X)$ is the distance to the boundary and $h_{r}$ is a one-parameter family of metrics on $\partial X$. We use $\{\alpha, \beta, \cdots\}$ as indices for objects on $\overline{X}$, $\{i,j, \cdots\}$ for objectors on $\partial X$ and $0$ for the $r$ variable. This means $\alpha$ corresponds to $(0, i)$ relative to the product identification induced by $\mathrm{exp}$. At $r=0$, the derivative $\partial^k_rh_r$ can be expressed in terms of curvature of $\bar{g}$. More explicitly,
\begin{equation*}
\begin{aligned}
&(h_r)|_{r=0} = \hat{g}, \\
&(\partial_rh_r)|_{r=0} =h', \quad h'_{ij}=-2L_{ij}, \\
&(\partial^2_rh_r)|_{r=0}=h'', \quad h''_{ij}=-2\bar{R}_{0i0j} +2L_{ik}L_j^{\ k},
\end{aligned}
\end{equation*}
where $L_{ij}$ is the second fundamental form on the boundary and $\bar{R}_{\alpha\beta\delta\gamma}$ is the Riemannian curvature of $\bar{g}$.

Let $\bar{R}$ and $\hat{R}$ be the scalar curvature of $\bar{g}$ and $\hat{g}$, respectively. Since $g_+$ has constant scalar curvature $-n(n-1)$, the conformal change of scalar curvature induces an equation for $\rho$:
\begin{equation}\label{eq.rho1}
-n(n-1) =-n(n-1)|\bar{\nabla}\rho|^2_{\bar{g}} +2(n-1)\rho \Delta_{\bar{g}}\rho +\rho^2 \bar{R},
\end{equation}
where $\Delta_{\bar{g}}=\bar{g}^{\alpha\beta}\bar{\nabla}_{\alpha}\bar{\nabla}_{\beta}$.
This implies that $\rho$ has an asymptotic expansion:
\begin{equation}\label{eq.asymrho}
\begin{aligned}
&\rho =r+c_2r^2 +c_3r^3+o(r^3),
\\
&c_2 = -\frac{1}{2(n-1)} H,
\\
&c_3 =\frac{1}{6(n-2)}(\hat{R}-|\mathring{L}|^2)-\frac{1}{6(n-1)} (\bar{R}+H^2),
\end{aligned}
\end{equation}
where $\mathring{L}_{ij}$ is the trace free part of $L_{ij}$ and $H$ is the mean curvature of the boundary.

Let $\bar{R}_{\alpha\beta}$ be the Ricci curvature tensor of $\bar{g}$ and $\bar{E}_{\alpha\beta}$ the trace free part of $\bar{R}_{\alpha\beta}$. Since $\mathrm{Ric}_{g_+}=-(n-1)g_+$, the conformal change of Ricci curvature implies that
\begin{equation}\label{eq.tracefreeRic}
\bar{E}=-(n-2)\rho^{-1} \left[\bar{\nabla}^2\rho - \frac{1}{n}(\Delta_{\bar{g}}\rho)\bar{g}\right].
\end{equation}
%

\vspace{0.2in}
\section{Proof of Theorem \ref{thm.2}}\label{sec3}

In this section, we prove Theorem \ref{thm.2}.

Let $(X^n, g_+)$ be a $C^{3,\alpha}$ conformally compact Poincar\'{e}-Einstein manifold satisfying the hypotheses of Theorem \ref{thm.Y}  with  first type Escobar-Yamabe compactification $(\overline{X}, \partial X, \bar{g})$. This means that $(\overline{X}, \bar{g})$ has constant scalar curvature $\bar{R}$ with totally geodesic boundary. Let $\hat{g}=\bar{g}|_{\partial X}$.
We further assume the equality in Theorem \ref{thm.GH} holds. Then  $\bar{g}$ is Einstein, i.e. the trace free part of Ricci curvature vanishes, and $\hat{g}$ has constant scalar curvature $\hat{R}$ on the boundary.

We will show that $(\overline{X}, \partial X, \bar{g})$ must be isometric to the standard hemisphere up to scalings. The idea here follows closely to Obata \cite{Ob1}, Cheeger-Colding \cite{ChC} and Wu-Ye \cite{WuYe}. Here we apply it in the case of a compact manifold with boundary with some extra conditions.

First, we apply the asymptotics in Section \ref{sec2} to the first type Escobar-Yamabe compactification $\bar{g}$. Since $H=0$, (\ref{eq.asymrho}) implies that in this case the boundary defining function has asymptotics
\begin{equation*}
\begin{gathered}
\rho =r+c_2r^2 +c_3r^3+o(r^3),
\end{gathered}
\end{equation*}
where
$$
c_2 = 0,
\quad
c_3 =\frac{1}{6(n-2)}\hat{R}-\frac{1}{6(n-1)} \bar{R}.
$$
Hence
\begin{equation}\label{eq.asym1}
\lim_{r \to 0} \frac{|\bar{\nabla} \rho|_{\bar g}^2 -1}{\rho}=4c_2=0.
\end{equation}
Rewrite equation (\ref{eq.rho1}) as
\begin{equation}\label{eq.deltarho}
\Delta_{\bar{g}} \rho =\frac{n}{2} \rho^{-1} \left(|\bar{\nabla} \rho|^2_{\bar{g}}-1\right)-\frac{1}{2(n-1)}\bar{R}\rho.
\end{equation}
Then (\ref{eq.asym1}) implies that on the boundary,
$$
|\bar{\nabla}\rho|_{\bar{g}}^2|_{\partial X}=1, \quad
\Delta_{\bar{g}} \rho|_{\partial X}=0.
$$
Since $\bar{E}=0$, by (\ref{eq.tracefreeRic}) we have
\begin{equation}\label{eq.rho2}
\bar{\nabla}^2\rho - \frac{1}{n}(\Delta_{\bar{g}}\rho)\bar{g}=0.
\end{equation}
For simplicity, we also denote
\begin{equation}\label{eq.z}
z=- \frac{1}{n}(\Delta_{\bar{g}}\rho)=-\frac{|\bar{\nabla}\rho|^2-1}{2\rho} +\frac{\bar{R}}{2n(n-1)} \rho .
\end{equation}
Hence $(\rho, z)$ satisfies the equation
\begin{equation}\label{eq.rho3}
\bar{\nabla}^2\rho +z\bar{g}=0.
\end{equation}
Since the boundary is totally geodesic, the Gauss-Codazzi equation
$
\bar{R}=\hat{R}+2\bar{R}_{00} + |L|^2-H^2
$
 gives that on the boundary
\begin{equation}\label{eq.scalar}
\hat{R}=\frac{n-2}{n}\bar{R}.
\end{equation}
Notice that up to now the sign of $\bar{R}$ (and hence $\hat{R}$) has not been determined.

Second, we consider a flow generated by the vector field $\bar{\nabla}\rho/|\bar{\nabla}\rho|^2$. Denote by $F(t,p)$ the flow lines of $\bar{\nabla}\rho/|\bar{\nabla}\rho|^2$ starting on $\partial X$, where $p\in \partial X$ is an initial point and $t$ is the time, i.e. $F(0,p)=p.$
Since $|\bar{\nabla}\rho|=1$ on $\partial X$, there exists some $T>0$ such that
$$
F:[0,T)\times \partial X\longrightarrow \overline{X}
$$
is a $C^{3,\alpha}$ diffeomorphism to its image and
$$|\bar{\nabla}\rho|(F(t,p))\neq 0, \quad \forall\ t\in [0,T).$$
This means $F(\{t\}\times\partial X)$ is a regular component of the level set of $\rho$.
We divide the rest part of the proof into the following claims. 
\begin{itemize}

\item[(a)] $\rho(F(t,p))=t$ for all $t\in [0,T)$.
\\
\textit{Proof.}
This is because $\rho(p)=0$ for $p \in \pa X$ and
$$
\frac{d}{dt}\rho(F(t,p))=\bar{\nabla}\rho \cdot \frac{\bar{\nabla}\rho}{|\bar{\nabla}\rho|^2}=1.
$$

\item[(b)] $|\bar{\nabla}\rho|(F(t,p))$ is independent of $p\in \partial X$ for all $t\in [0,T)$.
\\ 
\textit{Proof.}
This is because the equation  (\ref{eq.rho3}) is equivalent to
\begin{equation*}\label{eq.vector}
\bar{\nabla}_{V}\bar{\nabla} \rho +z V=0, \quad \forall\ V\in T\overline{X}.
\end{equation*}
If choosing $V$ to be tangent to $F(\{t\}\times\partial X)$, a regular level set of $\rho$, we have
$$
\bar{\nabla}_{V}|\bar{\nabla} \rho|^2 =-2z  V\cdot \bar{\nabla}\rho=0.
$$

\item[(c)] $[\triangle_{\bar{g}}\rho](F(t,p))$ and hence $z(F(t,p))$ are independent of $p\in \partial X$ for all $t\in [0,T)$.
\\
\textit{Proof.}
In view of equation (\ref{eq.deltarho}) and  by (a) and (b),  the function $[\triangle_{\bar{g}}\rho](F(t,p))$ is independent of $p$ for all $t\in [0,T)$.

\item[(d)]  $z(F(t,p)) =\bar{R}t/[n(n-1)], \ \forall\ t\in[0, T).$
\\
\textit{Proof.}
By \eqref{eq.z} and \eqref{eq.rho3}, we compute the $t$ derivative as follows
$$
\begin{aligned}
\frac{d}{dt}z(F(t,p)) & = \bar{\nabla}\left( -\frac{|\bar{\nabla}\rho|^2-1}{2\rho}+\frac{\bar{R}}{2n(n-1)}\rho
\right) \cdot  \frac{\bar{\nabla}\rho}{|\bar{\nabla}\rho|^2}
\\
& = \left(-\frac{\Delta_{\bar{g}}\rho}{n\rho} +\frac{|\bar{\nabla}\rho|^2-1}{2\rho^2}+\frac{\bar{R}}{2n(n-1)}
\right) \bar{\nabla}\rho\cdot  \frac{\bar{\nabla}\rho}{|\bar{\nabla}\rho|^2}
\\
&= \frac{\bar{R}}{n(n-1)},
\end{aligned}
$$
which is a constant. Since $z(F(0,p))= 0$ by \eqref{eq.asym1} and \eqref{eq.deltarho}, the claim then follows.

\item[(e)] The scalar curvature $\bar{R}>0$.
\\
\textit{Proof.}
Let $T_0$ be the maximum choice of $T$ such that $|\bar{\nabla}\rho|(F(t,p))\neq 0$ for all $t\in [0,T)$. $T_0$ is finite since $\rho$ is bounded on $\overline{X}$. Let $\Lambda$ be the limit set of the flow as $t\rightarrow T_0$, i.e.
\begin{equation}\label{limitset}
\Omega=F([0, T_0)\times \partial X), \quad
\Lambda = \overline{\Omega}-\Omega.
\end{equation}
It is obvious that $\Lambda\neq \emptyset$.

$\forall~q\in\Lambda$, we show that $q$ must be a critical point of $\rho$, i.e.  $|\bar{\nabla}\rho(q)|=0$. By definition, there exists $(t_k, p_k)\in [0,T_0)\times \partial X$ such that $F(t_k, p_k)\longrightarrow q$. By passing to a subsequence if necessary, we may further assume that $t_k$ increases to $T_0$ and $p_k$ converges to $p_0$. If $|\bar{\nabla}\rho(q)|\neq 0$,  then (b) implies that all the flow lines can be extended over $T_0$. This contradicts the maximum choice of $T_0$.

Next suppose $\bar{R}\leq 0$.

If $\bar{R}=0$, then
by (d), $z\equiv 0$ on the flow lines. By taking limit, we have $[\triangle_{\bar{g}}\rho](q)=0$. However, by (\ref{eq.deltarho}),
$[\triangle_{\bar{g}}\rho](q)=-n/[2\rho(q)] <0$, where we have used that $q$ is a critical point. Hence a contradiction.

If $\bar{R}<0$, we consider a global maximum point of $\rho$, denoted by $q_0$. Let
$\max_{\overline{X}}\rho =\rho_0=\rho(q_0)>0.$
Then by (\ref{eq.deltarho})
\begin{align} \notag
 -\frac{n}{2} -\frac{1}{2(n-1)}\bar{R}\rho_0^2=&[\rho\triangle_{\bar{g}}\rho](q_0)\leq 0.
\end{align}
This implies that
\begin{align}\notag
[\rho\triangle_{\bar{g}}\rho](q) &= -\frac{n}{2} -\frac{\bar{R}}{2(n-1)}\rho_0^2 +\frac{\bar{R}}{2(n-1)}(\rho_0^2-\rho^2(q))\leq 0.
\end{align}
However, in this case, (d) implies that $z<-C<0$ for all $t>T_0/2$. So $\rho\triangle_{\bar{g}} \rho$ is strictly positive at $q$ by taking limit. This yields a contradiction again.

\item[(f)] All critical points of $\rho$  are non-degenerate local maxima and $\overline{X}=\overline{\Omega}$. 
\\
\textit{Proof.}
By (d), (e) and (\ref{eq.rho3}), the Hessian of $\rho$ is negative definite at critical point $q\in \Lambda$, which is defined in (\ref{limitset}). So any point $q\in \Lambda$ is non-degenerate and local maximum. In particular, $\Lambda$ consists of  isolated critical points. For any $q\in \Lambda$, we can take a geodesic ball $B_r(q)$ with radius $r$ small enough such that $B_r(q)\cap \Lambda=\{q\}$. Since  $(B_r(q)-\{q\}) \cap \Omega\neq \emptyset$, we have $(B_r(q)-\{q\}) \subset \Omega$ , i.e. $q$ is an interior point of $\overline{\Omega}$. This implies $\overline{\Omega}$ is a smooth manifold with the only boundary $\partial X$. Since $\overline{X}$ is connected,  $\overline{X}=\overline{\Omega}$. 
This also shows that there are no other critical points except those  in $\Lambda$. 

\item[(g)] $\partial X$ is connected.
\\
\textit{Proof.}
Since $(\overline{X}, \partial X, \bar{g})$ is the first type Escobar-Yamabe compactification, $\bar{R}>0$ implies that $Y(\overline{X}, \partial X, [\bar{g}])>0$. Since the equality in Theorem \ref{thm.GH} holds in our case, we have $Y(\partial X, [\hat{g}])>0$ as well, which implies that $\partial X$ is connected by a theorem of Witten-Yau~\cite{WiYa}, see also ~\cite{CG}.

\item[(h)] $\overline{X}$ is diffeomorphic to a hemisphere and $\rho$ has only one critical point. 
\\
\textit{Proof.}
Since $\rho=0$ on $\partial X$, in view of claim (f) and standard Morse theory, we conclude that
\[
\overline{X}\cong (\partial X\times [0,1] )\bigsqcup_{\partial X\times \{1\}} (\sqcup_{i} D^{n}_i),
\]
i.e., $\overline{X}$ is resulted from attaching several copies of handles to $\partial X \times [0,1]$ along $\partial X\times\{1\}$. In this case, all handles are
$n$-disks. Therefore $\partial X$ is diffeomorphic to disjoint union of spheres. In view of (g), there exists only one disk, and the claim follows. Consequently, the only critical point of $\rho$ is where it attains its maximum on $\overline{X}$. Let us denote it by $q_0$.

\item[(i)] Level sets of $\rho$ are geodesic spheres of the unique critical point $q_0$.
\\
\textit{Proof.}
Combining claims (a) and (d), we can now rewrite (\ref{eq.rho3}) as
\begin{align}\label{eq.19}
\bar{\nabla}^2 \rho+ \frac{\bar{R}}{n(n-1)}\rho \bar{g}=0.
\end{align}
Up to a scaling, we may assume that $\bar{R}=n(n-1)$, thus (\ref{eq.19}) becomes to
\begin{align}\label{eq.20}
\bar{\nabla}^2 \rho+ \rho \bar{g}=0.
\end{align}
Taking any unit-speed geodesic $\gamma(s)$ starting from $q_0$, by (\ref{eq.20}) we have
\[
\frac{d^2}{ds^2} \rho(\gamma(s)) +  \rho(\gamma(s))=0,
 \quad \rho(\gamma(0))=T_0, \quad \frac{d}{ds}\rho(\gamma(s))\vline_{s=0}=0.
\]
It follows that
\[
\rho(\gamma(s))=T_0\cos (s),
\]
which means the geodesic $s$-sphere centered at $q_0$ is exactly the $T_0\cos(s)$-level set of $\rho$.
\end{itemize}

Finally, we piece all information together. According to the construction of the flow, and claims (h) and (i), it follows that
 $F: \partial X\times [0,T_0) \to \overline{X}\setminus \{q_0\}$ is a diffeomorphism. By (i), the flow lines coincide with geodesics starting at $q_0$ as they are both perpendicular to the level set of $\rho$. Thus equivalently, the exponential map at $q_0$, $\exp_{q_0}: B_{\frac{\pi}{2}}(0)\subset T_{q_0}\overline{X}\to X$ is a diffeomorphism. Using the coordinates associated with the exponential map at $q_0$, $\bar{g}$ can be expressed as
\[
\bar{g}= \ud s^2+ g_{s},
\] where $g_{s}$ is a smooth family of metrics on $S^{n-1}$, with $\displaystyle\lim_{s\to 0} s^{-2} g_s= g_{S^{n-1}}$.
Here $g_{S^{n-1}}$ is the standard round metric on $S^{n-1}$. Also by abuse of notation, $\rho(s)= T_0 \cos (s)$. As the variation of $g_s$ is the second fundamental form of geodesic $s$-sphere, we infer from (\ref{eq.20})
and $\rho(s)=T_0\cos(s)$ that
\[
\frac{d}{ds} g_s=2\cot(s) g_s.
\]
Hence $g_s=\sin^2 (s) g_0$. Then $(\overline{X}, \bar{g})$ is isometric to the standard hemisphere and $g_{+}= \rho^{-2}\bar{g}$ is the standard hyperbolic metric on $X$. This finishes the proof of Theorem \ref{thm.2}.

\vspace{0.2in}
\section{Proof of Theorem \ref{thm.1}} \label{sec4}
In this section we prove Theorem \ref{thm.1}.

Let $(X^n, g_+)$ be a $C^{3,\alpha}$ conformally compact Poincar\'{e}-Einstein manifold satisfying the hypotheses of Theorem \ref{thm.1} and hence those of Theorem \ref{thm.Q}.

We consider the second type Escobar-Yamabe compactification $(\overline{X}, \partial X, \bar{g})$. This means that $(\overline{X}, \bar{g})$ is scalar flat with constant mean curvature $H$ on the boundary. Such $\bar{g}$ may not be unique if exists. Here we just fix a choice of $\bar{g}=\rho^2g_+$ and let $\hat{g}=\bar{g}|_{\partial X}$.

Notice that in our setting, $(X, g_+)$ is $C^{3,\alpha}$ conformally compact. Hence by the boundary regularity theorem in~\cite{CDLS}, $(\overline{X}, \bar{g})$ has an umbilical boundary, i.e. the trace free part of $L_{ij}$, denoted by $\mathring{L}_{ij}$, vanishes on $\pa X$, and hence
$$
L_{ij}=\frac{H}{n-1}\hat{g}_{ij}.
$$

We apply the asymptotics in Section \ref{sec2} to the second type Escobar-Yamabe
compactification $\bar{g}$.
Since  $\bar{R}=0$ and  $\mathring{L}_{ij}=0$ for $\bar{g}$ now, (\ref{eq.asymrho}) implies that in this case the boundary defining function has an asymptotical expansion:
$$
\rho=r+c_2r^2 +c_3r^3+O(r^{3+\alpha}),
$$
where $r$ is the distance to the boundary and
\begin{align*}
c_2 =& -\frac{1}{2(n-1)}H,
\\
c_3 =&\frac{1}{6(n-2)}\hat{R}-\frac{1}{6(n-1)} H^2.
\end{align*}
Hence
\begin{equation}\label{eq.asym2}
\lim_{r \to 0} \frac{|\bar{\nabla} \rho|_{\bar g}^2 -1}{\rho}=4c_2=-\frac{2}{n-1}H.
\end{equation}

Because $\bar{R}=0$, we can rewrite equation (\ref{eq.rho1}) for $\rho$, which comes from the conformal change of scalar curvatures,  as
\begin{equation}\label{eq.deltarho2}
\Delta_{\bar{g}} \rho =\frac{n}{2} \rho^{-1} \left(|\nabla \rho|^2_{\bar{g}}-1\right).
\end{equation}
Recall  $\bar{E}$ is the trace-free part of Ricci  tensor of $\bar{g}$ and the conformal change of Ricci curvature gives
\begin{equation}\label{eq:trace-free_Ricci}
\bar{E} = -(n-2)\rho^{-1} \left[ \bar{\nabla}^2\rho - \frac{1}{n}(\Delta_{\bar{g}}\rho)\bar{g}\right].
\end{equation}

 Set $X_{\delta}:=\{x \in X| \ud_{\bar{g}}(x, \partial X)\geq \delta \}$ for small $\delta>0$. By \eqref{eq.deltarho2} and \eqref{eq:trace-free_Ricci}, a similar integration by parts argument as in \cite{GH} shows that
\begin{align}\label{eq:Obata1} 
\frac{2}{(n-2)^2} \int_{X_{\delta}} |\bar{E}|_{\bar g}^2 \rho \ud V_{\bar{g}} &=
\int_{\partial X_{\delta}} \rho^{-1}\left[  N(|\bar{\nabla}\rho|^2_{\bar{g}}) +\rho^{-1} (1-|\bar{\nabla}\rho|^2_{\bar{g}}) N(\rho)\right] \ud S_{\hat{g}},
\end{align}
where $N$ is the unit outward normal w.r.t. $\bar g$ on $\partial X_{\delta}$.

Clearly, $N=-\partial_r$.
Then a direct computation on the asymptotics of $\rho$ shows that
\begin{align*}
&\rho^{-1}=r^{-1}-c_2+(c_2^2-c_3)r+O(r^{1+\alpha}),
\\
&|\nabla \rho|_{\bar g}^2=1+4c_2 r+(4c_2^2+6c_3) r^2+O(r^{2+\alpha}).
\\
&\frac{|\nabla \rho|_{\bar g}^2 -1}{\rho}=4c_2+6c_3r+O(r^{1+\alpha}),
\\
&N(\rho)=-1-2c_2r-3c_3r^2+O(r^{2+\alpha}),
\\
&N(|\nabla \rho|_{\bar g}^2)=-4c_2-(8c_2^2+12c_3)r+O(r^{1+\alpha}).
\end{align*}
And hence
\begin{align*}
&\rho^{-1}\left[  N(|\nabla\rho|^2_{\bar{g}}) +\rho^{-1} (1-|\nabla\rho|^2_{\bar{g}}) N(\rho)\right]
=-6c_3+O(r^{\alpha})=\frac{1}{n-1} H^2 - \frac{1}{n-2}\hat{R}+O(r^{\alpha}).
\end{align*}
Letting $\delta \to 0$ in  \eqref{eq:Obata1}, we obtain
\begin{equation}\label{eq:Obata2}
\frac{2}{(n-2)^2} \int_X |\bar{E}|_{\bar g}^2 \rho  \ud V_{\bar{g}}=\int_{\partial X} \left(\frac{1}{n-1} H^2 - \frac{1}{n-2}\hat{R} \right) \ud S_{\hat{g}}.
\end{equation}

We are ready to prove the inequalities in Theorem \ref{thm.1} now:
\begin{itemize}
\item
When $n\geq 4$, it follows from \eqref{eq:Obata2} and Theorem \ref{thm.Q} that
\begin{equation}\label{est:Q_Y_n>3}
\begin{aligned}
Y(\pa X,[\hat g]) &\leq
\left( \int_{\pa X}\hat R \ud S_{\hat g} \right) \mathrm{Vol}(\pa X,\hat g)^{-\frac{n-3}{n-1}}\\
&=\left(\frac{n-2}{n-1}\int_{\partial X} H^2 \ud S_{\hat g}-\frac{2}{n-2} \int_X |\bar{E}|_{\bar g}^2 \rho \ud V_{\bar{g}}\right) \mathrm{Vol}(\pa X,\hat g)^{-\frac{n-3}{n-1}}\\\
&\leq \frac{n-2}{4(n-1)} Q(X,\pa X,[\bar g])^2.
\end{aligned}
\end{equation}
\item
When $n=3$, by Gauss-Bonnet theorem and \eqref{eq:Obata2} we have
\begin{equation}\label{est:Q_Y_n=3}
\begin{aligned}
4\pi \chi(\pa X)=\int_{\pa X}\hat R \ud S_{\hat g}=&\frac{1}{2}\int_{\partial X} H^2 \ud S_{\hat g}-2\int_X |\bar{E}|_{\bar g}^2 \rho \ud V_{\bar{g}}
\leq\frac{1}{8} Q(X,\pa X,[\bar g])^2.
\end{aligned}
\end{equation}
\end{itemize}

Next we consider the situation when the equality occurs. One direction is obvious. That is if $(X,g_+)$ is the standard hyperbolic manifold, then we can take a conformal compactification such that $(\overline{X}, \bar{g})$ is the Euclidean ball. A direct computation shows that the equality occurs in  \eqref{ineq:Q_Y_n>3} if $n>3$, or in \eqref{ineq:Q_Y_n=3} if $n=3$.  For the converse direction, we need to show that
when the equality occurs, the Escobar-Yamabe metric $\bar{g}$ must be the Euclidean metric and  $(\overline{X}, \bar{g})$ must be the Euclidean ball.

We observe from \eqref{est:Q_Y_n>3} and \eqref{est:Q_Y_n=3} that in both cases if the equality holds, then $\bar{E}=0$.
So the Ricci curvature of $\bar{g}$ vanishes, i.e. $\bar{R}_{\alpha\beta}=0$. Recall that  the the boundary is umbilic. These with  the Gauss-Codazzi equation
\begin{equation*}
\bar{R}=\hat{R}+2\bar{R}_{00}+ |L|_{\hat g}^2-H^2,
\end{equation*}
give that on the boundary
$$
\hat{R}=\frac{n-2}{n-1}H^2.
$$
Hence $\hat{R}\geq 0$ is also a constant. Moreover,  in this case, \eqref{eq:trace-free_Ricci} shows that
$$
\bar{\nabla}^2\rho = \frac{1}{n}(\Delta_{\bar{g}}\rho)\bar{g}.
$$
Thus we get
$$
\bar{\nabla}^{\beta}\bar{\nabla}_{\beta}\bar{\nabla}_{\alpha}\rho = \frac{1}{n} \bar{\nabla}_{\alpha}(\Delta_{\bar{g}}\rho).
$$
On the other hand,
$$
\bar{\nabla}^{\beta}\bar{\nabla}_{\beta}\bar{\nabla}_{\alpha}\rho
= \bar{\nabla}_{\alpha}\bar{\nabla}^{\beta}\bar{\nabla}_{\beta}\rho +\bar{R}_{\alpha\beta}\bar{\nabla}^{\beta}\rho
=\bar{\nabla}_{\alpha}(\Delta_{\bar{g}}\rho).
$$
Since $n\geq 3$,
$$
\bar{\nabla} (\Delta_{\bar{g}}\rho) =\frac{1}{n} \bar{\nabla}(\Delta_{\bar{g}}\rho)
\quad \Longrightarrow \quad  \bar{\nabla}(\Delta_{\bar{g}}\rho)=0.
$$
Therefore $ \Delta_{\bar{g}}\rho$ is a constant all over $\overline{X}$. By \eqref{eq.asym2} and \eqref{eq.deltarho2}  we have
\begin{equation}\label{eq.rho}
 \Delta_{\bar{g}}\rho =  (\Delta_{\bar{g}}\rho)|_{\partial X}= -\frac{n}{n-1} H.
\end{equation}
This implies that all over $\overline{X}$,
$$
\bar{\nabla}^2\rho = -\frac{1}{n-1}H\bar{g}.
$$

We claim that $H\neq 0$. Otherwise, $\Delta_{\bar{g}}\rho=0$ and $\rho|_{\partial X}=0$ implies that $\rho\equiv 0$ all over $\overline{X}$, which obviously can not happen. Now we can set $u=-(n-1)\rho/(nH)$. Then  $u$ satisfies
\begin{equation}\label{eq:u_PDE}
\begin{cases}
\Delta_{\bar g}u=1 &\hbox{~~in~~} X,\\
u=0 &\hbox{~~on~~} \pa X,
\end{cases}
\end{equation}
and $N(u)$ is constant on $\pa X$,
$$
N(u) = \frac{n-1}{nH}.
$$

Now the Reilly's formula in \cite{Reilly} together with \eqref{eq:u_PDE} gives
$$
\begin{aligned}
\frac{n-1}{n}\mathrm{Vol}(\overline{X}, \bar{g})
& =\frac{n-1}{n} \int_X (\Delta_{\bar{g}}u)^2 \ud V_{\bar{g}}
\\
&=\int_X\left[ (\Delta_{\bar{g}}u)^2 -|\bar{\nabla}^2u|^2_{\bar{g}} \right]\ud V_{\bar{g}}
\\
&=\int_{\partial X} H N(u)^2 \ud S_{\hat{g}}
\\
&=\left(\frac{n-1}{n}\right)^2\int_{\partial X} \frac{1}{H} \ud S_{\hat{g}}.
\end{aligned}
$$
Therefore we conclude that
$$
\int_{\partial X} \frac{n-1}{H} \ud S_{\hat{g}}=n\mathrm{Vol}(X, \bar{g}) .
$$
Hence $H>0$.
Then it follows from  \cite[Theorem 1]{Ros} that $(\overline{X}, \bar{g})$ is isometric to an Euclidean ball.

Up to a constant scaling, we can assume $H=n-1$. Then solving equation (\ref{eq.rho}) we can get $\rho=(1-r^2)/2$ and $(\overline{X}, \bar{g})$ is the unit ball in Euclidean space. And hence $(X, g_+=\rho^{-2}\bar{g})$ is the standard hyperbolic space $(\mathbb{H}^n,g_{\mathbb{H}})$. We finish the proof of Theorem \ref{thm.1}.

\vspace{0.2in}
\section{Proof of Corollary \ref{cor1}} \label{sec5}
In this section we prove Corollary \ref{cor1}.

Let  $(X,g_+)$ satisfy the same hypotheses as in Theoerem \ref{thm.Q}. Assume the conformal infinity of $(X, g_+)$ is conformally equivalent to the standard sphere $(S^{n-1}, g_{S^{n-1}})$.
Then
$$
\begin{aligned}
&Y(\pa X,[\hat g])=Y(S^{n-1},[g_{S^{n-1}}])\quad  &\textrm{if $n>3$};
\\
&\chi(\partial X)=\chi(S^2)& \quad \textrm{if $n=3$}.
\end{aligned}
$$
Combing this with the inequalities \eqref{ineq:Q_Y_n>3} and \eqref{ineq:Q_Y_n=3} in Theorem \ref{thm.1},  we obtain the following:
\begin{itemize}
 \item
 If $n>3$ then
\begin{equation}\label{eq.cor1}
\begin{aligned}
Y(S^{n-1},[g_{S^{n-1}}])=Y(\pa X,[\hat g])\leq& \frac{n-2}{4(n-1)} Q(X,\pa X,[\bar g])^2\\
\leq& \frac{n-2}{4(n-1)} Q(B^n,S^{n-1},[g_{\mathbb{R}^n}])^2=Y(S^{n-1},[g_{S^{n-1}}]).
\end{aligned}
\end{equation}
\item
If $n=3$, then
\begin{equation}\label{eq.cor2}
\begin{aligned}
32\pi\chi (S^2)=32\pi \chi (\pa X)\leq Q(X,\pa X,[\bar g])^2 \leq Q(B^3,S^{2},[g_{\mathbb{R}^3}])^2 = 32\pi\chi (S^2).
\end{aligned}
\end{equation}
\end{itemize}
These force all the inequalities in \eqref{eq.cor1} or \eqref{eq.cor2} to be equalities.
By Theorem \ref{thm.1}, this implies that $(X,g_+)$ is isometric to the hyperbolic space $(\mathbb{H}^n,g_{\mathbb{H}})$.
We finish the proof of  Corollary \ref{cor1}.

\end{document}